\documentclass[12pt,reqno]{amsart}
\usepackage{amsmath, amssymb, latexsym, amscd, amsthm,amsfonts,amstext}
\usepackage[top=3.5cm, bottom=3.5cm, left=3.5cm, right=3.5cm]{geometry} 
\usepackage[mathscr]{eucal}
\usepackage{xspace}
\usepackage{array, tabularx, longtable}
\usepackage{multicol,color}
\usepackage{indentfirst}
\usepackage{graphics}
\usepackage[colorlinks=true]{hyperref}
\theoremstyle{plain}
\newtheorem{theorem}{Theorem}

\newtheorem{lemma}[theorem]{Lemma}
\newtheorem{proposition}[theorem]{Proposition}
\theoremstyle{remark}

\newtheorem{remark}{Remark}

\theoremstyle{definition}
\newtheorem{definition}{Definition}

\begin{document}
\date{} 
\title [Non-integrated defect relations]{Non-integrated defect relations for the Gauss map of a complete minimal surface with finite total curvature in $\mathbb R^m$}

\author{Pham Hoang Ha}
\address{Department of Mathematics, Hanoi National University of Education, 136 XuanThuy str., Hanoi, Vietnam}
\email{ ha.ph@hnue.edu.vn}
\begin{abstract}
In this article, we give the non-integrated defect relations for the Gauss map of a complete minimal surface with finite total curvature in $\mathbb R^m.$ This is a continuation of previous work of Ha-Trao [J. Math. Anal. Appl., \textbf{430} (2015), 76-84.], which we extend here to targets of higher dimension.  
\end{abstract}
\maketitle
\section{Introduction}
In 1988, H. Fujimoto \cite{Fu2} proved Nirenberg's conjecture that if $M$ is a complete non-flat minimal surface in $\mathbb R^3,$ then its Gauss map can omit at most 4 points, and the bound is sharp. After that, he \cite{Fu4} also extended that result for minimal surfaces in $\mathbb R^m.$ He proved that the Gauss map of a non-flat complete minimal surface in $\mathbb R^m$ can omit at most $m(m+1)/2$ hyperplanes in $\mathbb P^{m-1}(\mathbb C)$ located in general position. He also gave an example to show that the number $m(m+1)/2$ is the best possible when $m$ is odd. Beside that, many mathematicians also studied the value distribution of the Gauss map of a minimal surface with finite total curvature and got many good results (see Fang \cite{Fa} and Ru \cite{Ru} for examples). On the other hand, Mo-Osserman \cite{MO} (1990), Mo \cite{M} (1994) and Ha-Phuong-Thoan \cite{HPT} recently showed the relations between the value distribution of the Gauss map and the total curvature of a complete minimal surface. Related to the value distributrion of the Gauss map of a complete minimal surface and the value distribution of the Gauss map of a complete minimal surface with finite total curvature, many results were given (see \cite{H}, \cite{DH},\cite{DHT}, \cite{KKM} and \cite{JR} for examples).\\
\indent On the other hand, Fujimoto \cite{Fu3, Fu4, Fu5} improved the previous results on the value distribution theory of the Gauss map of a complete minimal surface by introducing the modified defect relations for the Gauss map of a complete minimal surface which have analogy to the defect relations given by R. Nevanlinna in his value distribution theory. The author and Trao \cite{HT} recently improved the Fujimoto's results in the case the Gauss map of a complete minimal surface with finite total curvature in $\mathbb R^3, \mathbb R^4$ by studying the non-integrated defect relations for the Gauss map. In this article, we would like to be continuous to study the non-integrated defect relations for the Gauss map of a complete minimal surface with finite total curvature in $\mathbb R^m.$ These are the strict improvements of all previous results of Fujimoto on the modified defect relations for the Gauss map of a complete minimal surface with finite total curvature in $\mathbb R^m.$ Thus, they also are the improvements of previous results on ramifications for the Gauss map of a complete minimal surface with finite total curvature in $\mathbb R^m.$
\section{Statements of the main results}
 Let $M$ be an open Riemann surface and $f$ a nonconstant holomorphic map of $M$ into $\mathbb P^k(\mathbb C).$ Assume that $f$ has reduced representation $f=(f_0: \cdots : f_k).$ Set $||f|| = (|f_0|^2 +\cdots + |f_k|^2)^{1/2}$ and, for each a hyperplane $H : \overline{a}_{0}w_0+\cdots+ \overline{a}_{k}w_k = 0 $ in $\mathbb P^k(\mathbb C)$ with $|a_0|^2 +\cdots+ |a_k|^2 = 1,$ we definition $f(H) := \overline{a}_0f_0  +\cdots+\overline{a}_kf_k.$
\begin{definition} \label{d1}
We definition the $S-$defect of $H$ for $f$ by
$$\delta^{S}_{f}(H):= 1 - \inf \{\eta \geq 0; \eta \text{ satisfies condition $(*)_S$}\}.$$
Here, condition $(*)_S$ means that there exists a $ [-\infty, \infty)-$valued continuous subharmonic function $u\  (\not\equiv - \infty)$ on $M$ satisfying the following conditions:\\
(C1) $e^{u} \leq ||f||^{\eta},$ \\
(C2) for each $\xi \in f^{-1}(H) ,$ there exists the limit\\
$$\lim_{z \rightarrow \xi}(u(z) - \min (\nu_{f(H)}(\xi), k)\log|z-\xi|) \in [-\infty, \infty) ,$$
where $z$ is a holomorphic local coordinate around $\xi.$
\end{definition}
\begin{remark}
We always have that $\eta = 1$ satisfies condition $(*)_S$ when $u = \log |f(H)|.$
\end{remark}
\begin{definition}
\label{d2}
We definition the $H-$defect of $H$ for $f$ by
$$\delta^H_{f}(H):= 1 - \inf \{\eta \geq 0; \eta \text{ satisfies condition $(*)_H$}\}.$$
Here, condition $(*)_H$ means that there exists a $ [-\infty, \infty)-$valued continuous subharmonic function $u$ on $M$ which is harmonic on $M - f^{-1}(H)$ and satisfies the  conditions (C1) and (C2).
\end{definition}
\begin{definition} \label{d3}
We definition the $O-$defect of $H$ for $f$ by
$$\delta^O_{f}(H):= 1 - \inf \{\ \dfrac{1}{n}; \ \text{ $f(H)$ has no zero of order less than $n$}\}.$$
\end{definition}
\begin{remark}\label{rm1}
We always have $0\leq \delta^O_{f}(H ) \leq \delta^H_{f}(H)\leq \delta^S_{f}(H)\leq 1.$
\end{remark}
Moreover, Fujimoto \cite[page 672]{Fu1} also gave the reasons why he calls $\delta^S_f(H)$ the non-integrated defect by showing a relation between the non-integrated defect and the defect (as in Nevanlinna theory) of a nonconstant holomorphic map of $\Delta_R$ into $\mathbb P^k(\mathbb C).$
\begin{definition}
One says that $f$ is ramified over a hyperplane $H$ in $\mathbb P^{k}(\mathbb C)$ {\it with multiplicity at least} $e$ if all the zeros of the function $(f,H)$ have orders at least $e.$ If the image of $f$ omits $H,$ one will say that $f$ is {\it ramified over H with multiplicity }$\infty.$
\end{definition}
\begin{remark} \label{rm2}
If  $f$ is ramified over a hyperplane $H$ in $\mathbb P^{k}(\mathbb C)$ with multiplicity at least $n,$ then $\delta^S_f(H) \geq \delta^H_f(H) \geq \delta^O_f(H ) \geq 1 - \dfrac{1}{n}.$ In particular, if $f^{-1}(H) = \emptyset,$ then $\delta^O_f(H ) = 1.$ 
\end{remark}
In this article, we would like to study the $S-$defect relations for the Gauss maps of minimal surfaces with finite total curvature in  $\mathbb R^m$ which generalize the previous results of Ha-Trao in \cite{HT} to targets of higher dimension. In particular, we prove the following.\\
{\bf Main theorem.}\ {\it
Let $M$ be a non-flat complete minimal surface with finite total curvature in $\mathbb R^m$ and its Gauss map $G.$ Let $H_1,...,H_q$ be hyperplanes in $ \mathbb P^{m-1}(\mathbb C)$ located in $N$-subgeneral position $( q > 2N - k + 1, N \geq m-1).$ Assume that $G$ is $k-$non-degenerate (that is $G(M)$ is contained in a $k-$dimensional linear subspace in $ \mathbb P^{m-1}(\mathbb C)$, but none of lower dimension), $1\leq k \leq m-1,$ then
\begin{equation*}
 \sum_{j=1}^q\delta^S_{G}(H_j) \leq (k+1)(N-\dfrac{k}{2})+(N+1).
\end{equation*}  }
\begin{remark}
 For the case of the Gauss maps of minimal surfaces with finite total curvature, we can show that the Main theorem improved strictly Theorem 1.2 in \cite{Fu1}(by reducing the number $m^2$ to the number $m(m+1)/2$) and Theorem 2.8 in \cite{Fu4}(by changing the $H-$ defect relations to the $S-$ defect relations). 
\end{remark}
\begin{remark}
It is well known that the image of the (generalized) Gauss map $g: M \rightarrow  \mathbb P^{m-1}(\mathbb C)$ is contained in the hyperquadric $Q_{m-2}(\mathbb C) \subset \mathbb P^{m-1}(\mathbb C)$,
and that $Q_1(\mathbb C)$ is biholomorphic to $\mathbb P^1(\mathbb C)$ and that $Q_2(\mathbb C)$ is biholomorphic to $\mathbb P^1(\mathbb C)\times\mathbb P^1(\mathbb C) $. So the results Ha-Trao in (\cite{HT}) which only treat the cases $m=3$ and $m=4$ are better than a result which holds for any $m \geq 3$ can be if restricted to the special cases $m=3,4$. The easiest way to see the difference is to observe that $6$ lines in $\mathbb P^2(\mathbb C)$ in general position may have only $4$
points of intersection with the quadric $Q_1(\mathbb C) \subset \mathbb P^2(\mathbb C).$\\
\end{remark}
\section{Preliminaries and auxiliary lemmas}
In this section, we recall some auxiliary lemmas in \cite{Fu5, Fu6, Fu7}.\\
Let  $M$  be an open Riemann surface and $ds^2$ a pseudo-metric on $M$, namely, a metric on $M$ with isolated singularities which is locally written as $ds^2=\lambda ^2\left|dz\right|^2$ in terms of a nonnegative real-value function $\lambda$  with mild singularities and a holomorphic local coordinate $z$. We definition the divisor of  $ds^2$ by $\nu_{ds}:=\nu_{\lambda}$ for each local expression $ds^2=\lambda ^2 \left|dz\right|^2$, which is globally well-definitiond on $M$. We say that $ds^2$ is a continuous pseudo-metric if $\nu_{ds}\geq 0$ everywhere.
\begin{definition} (see \cite{Fu5})
We definition the Ricci form of $ds^2$ by 
$$\mbox{Ric}_{ds^2}:= -dd^c\log\lambda ^2$$
for each local expression $ds^2=\lambda ^2 \left|dz\right|^2.$ 
\end{definition}
In some cases, a $(1,1)-$form $\Omega$ on $M$ is regarded as a current on $M$ by defining $\Omega(\varphi):= \int_{M}\varphi\Omega$ for each $\varphi \in \mathcal{D},$ where $\mathcal{D}$ denotes the space of all $C^{\infty}$ differentiable functions on $M$ with compact supports.
\begin{definition} (see \cite{Fu5}) 
We say that a continuous pseudo-metric $ds^2$ has strictly negative curvature on $M$ if there is a positive constant $C$ such that 
$$- \mbox{Ric}_{ds^2}\geq C\cdot\Omega_{ds^2},$$
where $ \Omega_{ds^2}$ denotes the area form for $ds^2$, namely,
$$\Omega_{ds^2}:= \lambda ^2 (\sqrt{-1}/2)dz \wedge d\bar{z}$$
 for each local expression $ds^2=\lambda ^2\left| dz\right|^2$.
\end{definition}
As is well-known, if the universal covering surface of $M$ is biholomorphic with the unit disc in $\mathbb C,$ then $M$ has the complete conformal metric with constant curvature $-1$ which is called the Poincar\'{e} metric of $M$ and denoted by $d\sigma^2_M.$\\
\indent Let $f$ be a linearly non-degenerate holomorphic map of $M$ into $\mathbb P^k(\mathbb C).$ 
 Take a reduced representation $f = (f_0: \cdots : f_k)$. Then $F := (f_0, \cdots, f_k): M \rightarrow \mathbb C^{k+1} \setminus \{0\}$ is a holomorphic map with $\mathbb P(F) = f.$ 
 Consider the holomorphic map
\begin{equation*}
F_p=(F_p)_{z}:=F^{(0)}\wedge F^{(1)}\wedge\cdots\wedge F^{(p)}: M \longrightarrow \wedge^{p+1}\mathbb{C}^{k+1}
\end{equation*}
for $0\le p \le k,$  where  $F^{(0)} := F= (f_0,\cdots,f_k)$ and $F^{(l)}=(F^{(l)})_{z}:=(f_0^{(l)},\cdots,f_k^{(l)})$ for each $l=0,1,\cdots,k$, and where the $l$-th derivatives 
$f_i^{(l)}=(f_i^{(l)})_{z}$, $i=0,...,k$, are taken with respect to $z$.
(Here and for the rest of this paper the index $|_{z}$ means that the corresponding term
is definitiond by using differentiation with respect to the variable $z$, and in order to keep notations simple, we usually drop this index if no confusion is possible).
The norm of $F_p$ is given by
\begin{equation*}
\left| F_p\right|:= \bigg(\sum_{0\leq i_0<\cdots<i_p\leq k}\left|W(f_{i_0},\cdots,f_{i_p})\right|^2\bigg)^\frac{1}{2},
\end{equation*}
where $W(f_{i_0},\cdots,f_{i_p}) = W_z(f_{i_0},\cdots,f_{i_p})$ denotes the Wronskian of $f_{i_0},\cdots,f_{i_p}$ with respect to $z$.
\begin{proposition}(\cite[Proposition 2.1.6]{Fu7})\label{W}.\\
For two holomorphic local coordinates $z$ and $\xi$ and a holomorphic function 
$h : M \rightarrow \mathbb{C}$, the following holds :\\
a) $  W_{\xi}(f_0,\cdots, f_p)= W_z(f_0,\cdots, f_p) \cdot (\frac{dz}{d\xi})^{p(p+1)/2}$.\\
b) $W_{z}(hf_0,\cdots, hf_p)= W_z(f_0,\cdots, f_p) \cdot (h)^{p+1}. $ 
\end{proposition}
\begin{proposition}(\cite[Proposition 2.1.7]{Fu7})\label{W1}.\\
For holomorphic functions $f_0, \cdots , f_p : M \rightarrow \mathbb{C}$ the following conditions are equivalent:\\ 
(i)\ $ f_0, \cdots, f_p$ are linearly dependent over $\mathbb C.$\\
(ii)\ $W_z(f_0,\cdots,f_p) \equiv 0$ for some (or all) holomorphic local coordinate $z.$
\end{proposition}
We now take a hyperplane $H$ in $\mathbb P^k(\mathbb C)$ given by
\begin{equation*}
H:\overline{c}_0\omega_0+\cdots+\overline{c}_k\omega_k=0\,,
\end{equation*}
with $\sum_{i=0}^k|c_i|^2 = 1.$
We set 
\begin{equation*}
F_0(H) :=F(H):=\overline{c}_0f_0+\cdots+\overline{c}_kf_k
\end{equation*}
 and 
\begin{equation*}
\left| F_p(H)\right|=\left| (F_p)_{z}(H)\right|:= \bigg(\sum_{0\leq i_1< \cdots<i_p\leq k}\left|\sum_{l\not= i_1,...,i_p}\overline{c}_lW(f_{l},f_{i_1}, \cdots,f_{i_p})\right|^2\bigg)^\frac{1}{2},
\end{equation*}
for $1\leq p \leq k.$ We note that by using Proposition \ref{W}, $\left| (F_p)_{z}(H)\right|$ is multiplied by a factor $|\frac{dz}{d\xi}|^{p(p+1)/2}$ if we choose another holomorphic local coordinate $\xi$, and it is multiplied by $|h|^{p+1}$ if we choose another reduced representation $f=(hf_0:\cdots:hf_k)$ with a nowhere zero holomorphic function $h.$
Finally, for $0 \leq p \leq k$, set the $p$-th contact function of $f$ for $H$ to be $\phi_p(H):=\dfrac{|F_p(H)|^2}{|F_p|^2}=\dfrac{|(F_p)_{z}(H)|^2}{|(F_p)_{z}|^2}.$\\

We next consider $q$ hyperplanes $H_1,\cdots,H_q$ in $\mathbb{P}^{k}(\mathbb{C})$ given by
 $$H_j:\left\langle \omega ,A_j\right\rangle \equiv \overline{c}_{j0}\omega_0+\cdots+\overline{c}_{jk}\omega_k \quad(1\leq j\leq q)$$
where $A_j:=(c_{j0},\cdots,c_{jk})$ with $\sum_{i=0}^k|c_{ji}|^2 = 1.$

Assume now $N \geq k$ and $q \geq N+1$.
For $R\subseteq Q:=\left\{1,2,\cdots,q\right\},$ denote by $d(R)$ the dimension of the vector subspace of $\mathbb C^{k+1}$ generated by $\left\{A_j;j\in R \right\}$.\\
\indent The hyperplanes $H_1,\cdots,H_q$ are said to be in $N$-subgeneral position if $d(R)=k+1$ for all $R\subseteq Q$ with $\sharp (R)\geq N+1,$ where $\sharp (A)$ means the number of elements of a set $A.$ In the particular case $N=k$, these are said to be in general position.
\begin{theorem}(\cite[Theorem 2.4.11]{Fu7}) \label{N1}
\emph {\it For given hyperplanes $H_1,\cdots,H_q $ $( q > 2N - k + 1)$ in $\mathbb{P}^k(\mathbb{C})$ located in $N$-subgeneral position, there are some rational numbers $\omega(1),\cdots,\omega(q)$ and $\theta$ satisfying the following conditions:\\
 \indent (i) $0<\omega(j)\leq \theta\leq 1  \quad(1\leq j\leq q),$\\
 \indent (ii) $ \sum^{q}_{j=1}\omega(j)=k+1+\theta(q-2N+k-1),$\\
 \indent (iii) $\frac{k+1}{2N-k+1}\leq \theta \leq \frac{k+1}{N+1},$\\
 \indent (iv) If $R\subset Q$ and $0< \sharp(R)\leq n+1,$ then $\sum_{j \in R} \omega(j)\leq d(R).$}
 \end{theorem}
\noindent  Constants $\omega(j)\ (1 \leq  j  \leq q)$ and $\theta$ with the properties of Theorem~\ref{N1} are called Nochka weights and a Nochka constant for $H_1, \cdots, H_q$ respectively.

We need the three following results of Fujimoto combining the previously introduced concept of contact functions with Nochka weights:
\begin{theorem}(\cite[Theorem 2.5.3]{Fu7}) \label{PL1}
Let $H_1, \cdots, H_q$ be hyperplanes in $\mathbb P^k(\mathbb C)$ located in $N-$subgeneral position and let $\omega(j)$ $(1\leq j \leq q)$ and $\theta$ be Nochka weights and a Nochka constant for these hyperplanes. For every $\epsilon > 0$ there exist some positive numbers $\delta(>1)$ and $C,$ depending only on $\epsilon$ and $H_j\,$, $1\leq j \leq q,$ such that
\begin{align}\label{eq:T1}
&dd^c\log \dfrac{\Pi_{p=0}^{k-1}|F_p|^{2\epsilon}}{\Pi_{1\leq j\leq q, 0\leq p\leq k-1}\log^{2\omega(j)}(\delta/\phi_p(H_j))}\nonumber \\
&\geq C\bigg( \dfrac{|F_0|^{2\theta(q-2N+k-1)}|F_k|^2}{\Pi_{j=1}^q(|F(H_j)|^2\Pi_{p=0}^{k-1}\log^2(\delta/\phi_p(H_j)))^{\omega(j)}}\bigg)^{\frac{2}{k(k+1)}}dd^c|z|^2.
\end{align}
\end{theorem}
\begin{proposition}(\cite[Proposition 2.5.7]{Fu7})\label{P} \
 Set $\sigma_p=p(p+1)/2$ for $0 \leq p \leq k$ and $\tau_k = \sum_{p=0}^k\sigma_p.$ Then,
\begin{align}
dd^c\log(|F_0|^2|F_1|^2\cdots|F_{k-1}|^2)\geq \dfrac{\tau_k}{\sigma_k}\bigg(\dfrac{|F_0|^2|F_1|^2\cdots|F_{k}|^2}{|F_0|^{2\sigma_{k+1}}}\bigg)^{1/\tau_k}dd^c|z|^2.
\end{align}
\end{proposition}
\begin{proposition} (\cite[Lemma 3.2.13]{Fu7})
\label{P7}
 Let $f$ be a non-degenerate holomorphic map of a domain in  $\mathbb{C}$ into $\mathbb{P}^k(\mathbb{C})$ with reduced representation $f=(f_0:\cdots :f_k)$ and let $H_1,\cdots,H_q$ be hyperplanes located in $N$-subgeneral position $( q > 2N - k + 1)$ with Nochka weights $\omega(1),\cdots,\omega(q)$ respectively. Then,
$$\nu_\phi + \sum^{q}_{j=1}\omega(j)\cdot\min (\nu_{(f, H_j)}, k ) \geq 0,$$
where $\phi = \dfrac{|F_k|}{\Pi_{j=1}^q \mid F(H_j)\mid^{\omega (j)}}.$
\end{proposition}
\begin{lemma}(Generalized Schwarz's Lemma \cite{Ah}) \label{L3}
Let $v$ be a non-negative real-valued continuous subharmonic function on $\Delta_R.$ If $v$ satisfies the inequality  $\Delta\log v \geq v^2$ in the sense of distribution, then
$$v(z) \leq \dfrac{2R}{R^2 - |z|^2}.$$
\end{lemma}
\section{ The proof of the Main Theorem }
\begin{proof}
\indent For the convenience of the reader, we first recall some notations on the Gauss map of minimal surfaces in $\mathbb R^m.$
Let $M$ be a complete immersed minimal surface in $\mathbb R^m.$ Take an immersion $x
=(x_0,...,x_{m-1}) : M \rightarrow \mathbb R^m.$ Then $M$ has the structure of a Riemann surface and any local isothermal coordinate $(x, y)$ of $M$ gives a local holomorphic coordinate $z=x+ \sqrt{-1}y$. The generalized Gauss map of $x$ is definitiond to be 
\begin{equation*}\label{eq:}
g : M \rightarrow \mathbb P^{m-1}(\mathbb C), g = \mathbb P(\dfrac{\partial x}{\partial z})=(\dfrac{\partial x_0}{\partial z}:\cdots:\dfrac{\partial x_{m-1}}{\partial z}).
\end{equation*}
Since $x:M \rightarrow \mathbb R^m$ is immersed, $$G=G_{z}:= 
(g_0,...,g_{m-1}) = ((g_0)_{z},...,(g_{m-1})_{z}) =(\dfrac{\partial x_0}{\partial z},\cdots,\dfrac{\partial x_{m-1}}{\partial z})$$ is a (local) reduced representation of $g$, and
since for another local holomorphic coordinate $\xi$ on $M$ we have $G_{\xi} = G_{z}\cdot (\dfrac{dz}{d\xi})$, $g$ is well definitiond (independently of the (local) holomorphic coordinate). 
Moreover, if $ds^2$ is the metric on $M$ induced by the standard metric on $\mathbb R^m$, we have 
\begin{equation}\label{eq:1}
ds^2 = 2|G_{z}|^2|dz|^2 .
\end{equation}
Finally since $M$ is minimal,  $g$ is a holomorphic map. 

Since by hypothesis of the Main theorem, $g$ is $k$-non-degenerate $(1 \leq k \leq m-1)$
 without loss of
generality, we may assume that $g(M) \subset \mathbb P^k(\mathbb C);$ then
\begin{equation*}\label{eq:}
g : M \rightarrow \mathbb P^{k}(\mathbb C), g = \mathbb P(\dfrac{\partial x}{\partial z})=(\dfrac{\partial x_0}{\partial z}:\cdots:\dfrac{\partial x_{k}}{\partial z})
\end{equation*}
is linearly non-degenerate in $\mathbb P^k(\mathbb C)$ (so in particular $g$ is not constant)  and the other facts mentioned above still hold.

Let $H_j(j = 1,...,q)$ be $q (\geq N+1)$ hyperplanes in $\mathbb P^{m-1}(\mathbb C)$ in $N$-subgeneral position $(N\geq m-1 \geq k)$. Then  $H_j\cap \mathbb P^{k}(\mathbb C)(j = 1,...,q)$ are $q$ hyperplanes in $\mathbb P^{k}(\mathbb C)$ in $N$-subgeneral position. Let each $H_j\cap \mathbb P^{k}(\mathbb C)$ be represented  as
\begin{equation*}
H_j\cap \mathbb P^{k}(\mathbb C) :\overline{c}_{j0}\omega_0+\cdots+\overline{c}_{j{k}}\omega_{k}=0
\end{equation*}
with $\sum_{i=0}^k|c_{ji}|^2 = 1.$\\
Set
\begin{equation*}
G(H_j) = G_{z}(H_j):= \overline{c}_{j0}g_0+\cdots+\overline{c}_{jk}g_{k}.
\end{equation*}

We will now, for each contact function $\phi_p(H_j)$ for each of our hyperplanes $H_j$, choose one of the components of the numerator $|((G_z)_p)_z(H_j)|$ which is not identically zero: More precisely, for each $j,p\ (1\leq j \leq q, 1\leq p \leq k),$ we can choose $i_1,\cdots,i_p$ with $0\leq i_1<\cdots< i_p \leq k$ such that
\begin{equation*}
\psi(G)_{jp}=(\psi (G_z)_{jp})_{z}:=\sum_{l\neq i_1,..,i_p}\overline{c}_{jl}W_{z}(g_l,g_{i_1},\cdots,g_{i_p})\not\equiv 0,
\end{equation*}
(indeed, otherwise, we have 
 $\sum_{l\neq i_1,..,i_p}\overline{c}_{jl}W(g_l,g_{i_1},\cdots,g_{i_p})\equiv 0$ for all $i_1, ..., i_p$,
so $W(\sum_{l\neq i_1,..,i_p}\overline{c}_{jl}g_l,g_{i_1},\cdots,g_{i_p})\equiv 0$ for all $i_1, ..., i_p$, which contradicts the non-degeneracy of $g$ in $\mathbb P^k(\mathbb C).$
Alternatively we simply can observe that in our situation none of the contact functions vanishes identically).

Now we prove the Main theorem in four steps:\\
{\bf Step 1:} We may assume that
\begin{equation}\label{eq:2}
\sum_{j=1}^q\delta_{g}^S(H_j) > (k+1)(N-\dfrac{k}{2})+(N+1),  
\end{equation} 
otherwise our Main theorem is already proved. By definition, there exist constants $\eta_j \geq 0 (1\leq j \leq q)$ such that $q- \sum_{j=1}^q\eta_j > (k+1)(N-\dfrac{k}{2})+(N+1)$ and continuous subharmonic functions $u_j (1\leq j \leq q)$ on $M$ satisfies conditions (C1) and (C2). Thus
\begin{equation}\label{eq:2'}
\sum_{j=1}^{q}(1- \eta_j)-2N+k-1 > \dfrac{(2N-k +1)k}{2}>0,
\end{equation}
and this implies in particular
\begin{equation} \label{ass1''}
q>2N-k+1 \geq N+1 \geq k+1 .
\end{equation}
 By Theorem \ref{N1}, we have 
\begin{equation*}
 (q - 2N + k - 1)\theta = \sum_{j=1}^q \omega (j) - k - 1,\  \theta \geq \omega(j) > 0 
\text{ and } \theta \geq \dfrac{k + 1}{2N - k + 1},
\end{equation*}
so
\begin{align*}
2\bigg (\sum_{j=1}^q \omega(j)(1-\eta_j) - k - 1\bigg )&= 2(\sum_{j=1}^{q}\omega(j)- k - 1) - 2\sum_{j=1}^{q}\omega(j)\eta_j\\
&= 2(q-2N + k - 1)\theta - 2\sum_{j=1}^{q}\omega(j)\eta_j\\
&\geq 2(q-2N + k - 1)\theta - 2\sum_{j=1}^{q}\theta\eta_j\\
&= 2\theta\bigg (\sum_{j=1}^{q}(1- \eta_j)-2N+k-1\bigg )\\
&\geq 2\dfrac{(k + 1)\bigg (\sum_{j=1}^{q}(1- \eta_j)-2N+k-1\bigg )}{2N-k+1}.
\end{align*}
Thus, we now can  conclude with (\ref{eq:2'}) that 
\begin{align}\label{eq:3}
&2\bigg (\sum_{j=1}^q \omega (j)(1-\eta_j) - k - 1\bigg ) > k(k+1) \nonumber \\
&\Rightarrow \sum_{j=1}^q \omega (j)(1-\eta_j) - k - 1 - \dfrac{k(k+1)}{2} >0.
\end{align}
We set $\gamma := \sum_{j=1}^q \omega (j)(1-\eta_j) - k - 1.$\\
Then, by (\ref{eq:3}), we get
$\gamma > \sigma_{k} \Rightarrow \dfrac{\gamma}{\sigma_k} > 1.$\\
So we can choose a positive real number $\epsilon$ such that $\dfrac{\gamma-\epsilon\sigma_{k+1}}{\sigma_k + \epsilon\tau_k} > 1.$\\
{\bf Step 2:}\\
We set
\begin{equation*}
\lambda_z := \bigg( \dfrac{|G_z|^{\gamma -\epsilon\sigma_{k+1}}e^{\sum_{j=1}^q\omega(j)u_j}.|(G_k)_z|.\prod_{p=0}^{k}|(G_p)_z|^{\epsilon}}{\prod_{j=1}^{q}(|G(H_j)|\Pi_{p=0}^{k-1}\log(\delta/\phi_p(H_j)))^{\omega (j)}}\bigg)^{\frac{1}{\sigma_k +\epsilon\tau_k}},
\end{equation*}
and definition the pseudometric $d\tau_z^2 := \lambda_z^2|dz|^2.$ Using Proposition \ref{W}, we can see that 
\begin{align*}
d\tau_{\xi} &:= \bigg( \dfrac{|G_{\xi}|^{\gamma -\epsilon\sigma_{k+1}}e^{\sum_{j=1}^q\omega(j)u_j}.|(G_k)_{\xi}|.\prod_{p=0}^{k}|(G_p)_{\xi}|^{\epsilon}}{\prod_{j=1}^{q}(|G(H_j)|\Pi_{p=0}^{k-1}\log(\delta/\phi_p(H_j)))^{\omega (j)}}\bigg)^{\frac{1}{\sigma_k +\epsilon\tau_k}}|d\xi|\\
&= \bigg( \dfrac{|G_z|^{\gamma -\epsilon\sigma_{k+1}}e^{\sum_{j=1}^q\omega(j)u_j}.|(G_k)_{z}||\frac{dz}{d\xi}|^{\sigma_k}.\prod_{p=0}^{k}|(G_p)_{z}|^{\epsilon}.|\frac{dz}{d\xi}|^{\sum_{p=0}^{k}\epsilon\frac{p(p+1)}{2}}}{\prod_{j=1}^{q}(|G(H_j)|\Pi_{p=0}^{k-1}\log(\delta/\phi_p(H_j)))^{\omega (j)}}\bigg)^{\frac{1}{\sigma_k +\epsilon\tau_k}}|\dfrac{d\xi}{dz}|.|dz|\\
&= \bigg( \dfrac{|G_z|^{\gamma -\epsilon\sigma_{k+1}}e^{\sum_{j=1}^q\omega(j)u_j}.|(G_k)_{z}|.\prod_{p=0}^{k}|(G_p)_{z}|^{\epsilon}.|\frac{dz}{d\xi}|^{\sigma_k +\epsilon\tau_k}}{\prod_{j=1}^{q}(|G(H_j)|\Pi_{p=0}^{k-1}\log(\delta/\phi_p(H_j)))^{\omega (j)}}\bigg)^{\frac{1}{\sigma_k +\epsilon\tau_k}}|\dfrac{d\xi}{dz}|.|dz|\\
&=d\tau_{z}.
\end{align*}
Thus $d\tau_z^2$ is independent of the choice of the local coordinate $z.$ We will denote $d\tau_z^2$ by $d\tau^2$ for convenience. So $d\tau^2$ is well-definitiond on $M$ and 
\begin{equation*}
d\tau^2 =\bigg( \dfrac{|G|^{\gamma -\epsilon\sigma_{k+1}}e^{\sum_{j=1}^q\omega(j)u_j}.|G_k|.\prod_{p=0}^{k}|G_p|^{\epsilon}}{\prod_{j=1}^{q}(|G(H_j)|\Pi_{p=0}^{k-1}\log(\delta/\phi_p(H_j)))^{\omega (j)}}\bigg)^{\frac{2}{\sigma_k +\epsilon\tau_k}}|dz|^2:= \lambda^2|dz|^2.
\end{equation*}
{\bf Step 3:}\ We will show that $d\tau^2$ is continuous and has strictly negative curvature on $M$ in this step. 

Indeed, it is easy to see that $d\tau$ is continuous at every point $z_0$ with $\Pi_{j=1}^qG(H_j)(z_0) \not=0.$ Now we take a point $z_0 $ such that $\Pi_{j=1}^qG(H_j)(z_0) =0.$ Hence, it follows from (C2) that we get 
\begin{align*}
\lim_{z \to z_0} e^{u_j(z)}|z-z_0|^{-\min(\nu_{G(H_j)}(z_0), k)} =\lim_{z \to z_0} e^{u_j(z) - \min(\nu_{G(H_j)}(z_0), k)\log|z-z_0|} < \infty, \forall \  1\leq j \leq q.
\end{align*}
Thus, combining this with Proposition \ref{P7}, we get
\begin{align*}
\nu_{d\tau}(z_0)& \geq \dfrac{1}{\sigma_k + \epsilon\tau_k}\bigg( \nu_{G_k}(z_0)-\sum_{j=1}^q\omega(j)\nu_{G(H_j)}(z_0)+\sum_{j=1}^q\omega(j)\min\{\nu_{G(H_j)}(z_0), k \}\\
& +\sum_{j=1}^q\omega(j)\nu_{e^{u_j(z)}|z-z_0|^{-\min(\nu_{G(H_j)}(z_0), k)}}(z_0) \bigg) \\ 
&\geq 0.
\end{align*}
This concludes the proof that $d\tau$ is continuous on $M.$

On the other hand, by using Proposition \ref{P}, Theorem \ref{PL1} and noting that $dd^c\log|G_k| =0, dd^c\log e^{u_j} \geq 0 (1\leq j\leq q)$, we have
\begin{align*}
dd^c\log\lambda&\geq \dfrac{\gamma - \epsilon\sigma_{k+1}}{\sigma_{k}+\epsilon\tau_k}dd^c\log|G|+\dfrac{\epsilon}{4(\sigma_{k}+\epsilon\tau_k)}dd^c\log(|G_0|^2\cdots|G_{k-1}|^2)\\ 
& + \dfrac{1}{2(\sigma_k +\epsilon\tau_k)}dd^c\log\dfrac{\prod_{p=0}^{k-1}|G_p|^{2(\frac{\epsilon}{2})}}{\prod_{j=1}^{q}\Pi_{p=0}^{k-1}\log^{2\omega (j)}(\delta/\phi_p(H_j))}\\
&\geq \dfrac{\epsilon}{4(\sigma_{k}+\epsilon\tau_k)}\dfrac{\tau_k}{\sigma_k}\bigg(\dfrac{|G_0|^2|G_1|^2\cdots|G_{k}|^2}{|G_0|^{2\sigma_{k+1}}}\bigg)^{1/\tau_k}dd^c|z|^2\\ 
& + C_0\bigg( \dfrac{|G_0|^{2\theta(q-2N+k-1)}|G_k|^2}{\Pi_{j=1}^q(|G(H_j)|^2\Pi_{p=0}^{k-1}\log^2(\delta/\phi_p(H_j)))^{\omega(j)}}\bigg)^{\frac{2}{k(k+1)}}dd^c|z|^2\\
&\geq \min \{\dfrac{1}{4\sigma_k(\sigma_{k}+\epsilon\tau_k)}, \dfrac{C_0}{\sigma_k} \}\bigg(\epsilon\tau_k\bigg(\dfrac{|G_0|^2|G_1|^2\cdots|G_{k}|^2}{|G_0|^{2\sigma_{k+1}}}\bigg)^{1/\tau_k}\\
&+ \sigma_k\bigg( \dfrac{|G_0|^{2\theta(q-2N+k-1)}|G_k|^2}{\Pi_{j=1}^q(|G(H_j)|^2\Pi_{p=0}^{k-1}\log^2(\delta/\phi_p(H_j)))^{\omega(j)}}\bigg)^{\frac{1}{\sigma_k}}\bigg) dd^c|z|^2
\end{align*}
where $C_0$ is the positive constant. So,  by using the basic inequality 
$$ \alpha A + \beta B \geq (\alpha + \beta) A^{\frac{\alpha}{\alpha+ \beta}}B^{\frac{\beta}{\alpha + \beta}} \text{ for all }\alpha, \beta, A, B > 0  ,$$ we can find a positive constant $C_1$ satisfing the following\\
\begin{align*}
dd^c\log\lambda &\geq C_1\bigg( \dfrac{|G|^{\theta (q-2N+k-1) -\epsilon\sigma_{k+1}}.|G_k|.\prod_{p=0}^{k}|G_p|^{\epsilon}}{\prod_{j=1}^{q}(|G(H_j)|\cdot \Pi_{p=0}^{k-1}\log(\delta/\phi_p(H_j)))^{\omega (j)}}\bigg)^{\frac{2}{\sigma_k +\epsilon\tau_k}}dd^c|z|^2\\
&= C_1\bigg( \dfrac{|G|^{\sum_{j=1}^q \omega (j) - k-1 -\epsilon\sigma_{k+1}}.|G_k|.\prod_{p=0}^{k}|G_p|^{\epsilon}}{\prod_{j=1}^{q}(|G(H_j)|\cdot \Pi_{p=0}^{k-1}\log(\delta/\phi_p(H_j)))^{\omega (j)}}\bigg)^{\frac{2}{\sigma_k +\epsilon\tau_k}}dd^c|z|^2 \  \text{ (by Theorem \ref{N1}) }\\
&= C_1\bigg( \dfrac{|G|^{\gamma -\epsilon\sigma_{k+1}}.e^{\sum_{j=1}^q\omega(j)u_j}.|G_k|.\prod_{p=0}^{k}|G_p|^{\epsilon}.\prod_{j=1}^{q}\bigg(\dfrac{|G|^{\eta_j}}{e^{u_j}}\bigg)^{\omega (j)}}{\prod_{j=1}^{q}(|G(H_j)|\cdot \Pi_{p=0}^{k-1}\log(\delta/\phi_p(H_j)))^{\omega (j)}}\bigg)^{\frac{2}{\sigma_k +\epsilon\tau_k}}dd^c|z|^2.
\end{align*}
On the other hand, 
\begin{equation*}
\bigg(\dfrac{|G|^{\eta_j}}{e^{u_j}}\bigg)^{\omega (j)} \geq 1 \text{ for all } j = 1, ..., q,
\end{equation*}
so we get
\begin{align*}
dd^c\log\lambda &\geq C_1\bigg( \dfrac{|G|^{\gamma -\epsilon\sigma_{k+1}}.e^{\sum_{j=1}^q\omega(j)u_j}.|G_k|.\prod_{p=0}^{k}|G_p|^{\epsilon}}{\prod_{j=1}^{q}(|G(H_j)|\cdot \Pi_{p=0}^{k-1}\log(\delta/\phi_p(H_j)))^{\omega (j)}}\bigg)^{\frac{2}{\sigma_k +\epsilon\tau_k}}dd^c|z|^2\\
&=C_1\lambda^2dd^c|z|^2.
\end{align*}
This concludes the proof that $d\tau^2$ has strictly negative curvature on $M.$\\
{\bf Step 4:}\ Set $\Omega_g = dd^c\log||G||^2.$ \\
\indent Now, by using the classical Nevanlinna theory \cite[Theorem 3.3.15]{Fu7} and remarks of Fujimoto \cite[Proposition 4.7]{Fu1} we get $\sum_{j=1}^q\delta^S_g(H_j) \leq 2N-k+1$ if the universal covering surface of $M$  is biholomorphic to $\mathbb C.$ This is a contradiction with (\ref{eq:2}). Thus, we only need consider the case that the universal covering surface of $M$  is biholomorphic to the unit disc. By Lemma \ref{L3}, there exists a positive constant $C_0$ such that
 $$d\tau^2\leq C_0d\sigma^2_M,$$ where $d\sigma^2_M$ denotes the Poincar\'e metric on $M$.\\
Now, it follows from the assumption $M$ has finite total curvature in $\mathbb R^m$ that $M$ is bi-holomorphic with a compact Riemann surface $\overline{M}$ with finitely many points $a_l $'s removed. For each $a_l$, we take a neighborhood $U_l$ of $a_l$ which is biholomorphic to $\Delta^*=\left\{z; 0<|z| <1\right\},$ where $z(a_l)=0$. The Poincar\'e metric on domain $\Delta^*$ is given by
$$d\sigma^2_{\Delta^*}=\frac{4|dz|^2}{|z|^2\log^2|z|^2}.$$
By using the distance decreasing property of the Poincar\'e metric, we have
$$d\tau^2\leq C_l\frac{|dz|^2}{|z|^2\log^2|z|^2}$$
with some $C_l>0$ . This implies that, for a neighborhood $U^*_l$  of $a_l$ which is relatively compact in $U_l,$ we have
$$\int_{U^*_l}\Omega_{d\tau^2}<+\infty.$$
Since $\overline{M}$ is compact, we have
\begin{equation}\label{eq:1.1}
\int_M\Omega_{d\tau^2} \leq \int_{\overline{M}-\cup_l U^*_l}\Omega_{d\tau^2}+\sum_l\int_{U^*_l}\Omega_{d\tau^2}<+\infty.
\end{equation}
\indent On the other hand, we have 
\begin{align*}
dd^c\log\lambda&\geq (\dfrac{\gamma - \epsilon\sigma_{k+1}}{\sigma_{k}+\epsilon\tau_k})dd^c\log|G|+\dfrac{1}{(\sigma_k +\epsilon\tau_k)}dd^c\log\dfrac{\prod_{p=0}^{k-1}|G_p|^{\epsilon}}{\prod_{j=1}^{q}\Pi_{p=0}^{k-1}\log^{2\omega (j)}(\delta/\phi_p(H_j))}\\
&\geq (\dfrac{\gamma - \epsilon\sigma_{k+1}}{\sigma_{k}+\epsilon\tau_k})dd^c\log|G|, \text{ by Theorem \ref{PL1}.}
\end{align*}
 This implies
$$dd^c\log \lambda^2 \geq (\dfrac{\gamma - \epsilon\sigma_{k+1}}{\sigma_{k}+\epsilon\tau_k})\Omega_g.$$
Thus, we now can find a subharmonic function $v$ such that 
\begin{align*}
\lambda^2|dz|^2& = e^{v}||G||^{2(\dfrac{\gamma - \epsilon\sigma_{k+1}}{\sigma_{k}+\epsilon\tau_k})}|dz|^2\\
&=e^{v + (\dfrac{\gamma - \epsilon\sigma_{k+1}}{\sigma_{k}+\epsilon\tau_k} - 1)\log ||G||^2}||G||^2|dz|^2\\ 
&= e^w ds^2.
\end{align*}
So $d\tau^2 = e^wds^2,$ where $w$ is a subharmonic function. Here, we can apply the result of  Yau in \cite{Yau} to see
$$\int_M e^w\Omega_{ds^2} = + \infty,$$
because of the completeness of $M$ with respect to the metric $ds^2.$ This contradicts the assertion (\ref{eq:1.1}). The proof of the Main theorem is completed.
\end{proof}
{\bf Acknowledgements.} 
This work is supported by a grant of Hanoi National University of Education. The author also is grateful to Professors Do Duc Thai and Gerd Dethloff for many stimulating discussions concerning this material. 

\end{document}